\documentclass[amscd]{amsart}

\newtheorem{theorem}{Theorem}[section]
\newtheorem{lemma}[theorem]{Lemma}
\newtheorem{corollary}[theorem]{Corollary}

\usepackage{amsfonts}
\usepackage{amssymb}
\usepackage{amsmath}
\usepackage{eucal}
\usepackage{framed}
\usepackage[dvips]{graphicx}
\usepackage{epic}
\usepackage{graphics}
\usepackage{amsthm}
\usepackage{fancybox}
\usepackage{color} 
\usepackage{ulem} 
\usepackage{eepic}

\begin{document}

\begin{flushright}
  {\tiny --[GGJAC2022-10-26-Some-Comments-withFig-rf.tex)]--}
\end{flushright}

\bigskip
  
\title{Some Comments around The Examples against\\
 The Generalized Jacobian Conjecture}

\author{Susumu ODA}
\maketitle

{\small
 \begin{center}
Professor Emeritus, Kochi University\\   
{\small  345  Hiradani, Taki-cho, Taki-gun, Mie 519-2178, JAPAN\\
   ssmoda@ma.mctv.ne.jp}
\end{center}
}


\begin{abstract} 
 We have studied a faded problem, the Jacobian Conjecture ~: 

\noindent
  {\sf The Jacobian Conjecture $(JC_n)$}~:
 If  $f_1, \cdots, f_n$  are elements in a polynomial ring  $k[X_1, \cdots, X_n]$  over a field $k$ of characteristic $0$ such that the Jacobian $\det(\partial f_i/ \partial X_j) $ is a nonzero constant, then  $k[f_1, \cdots, f_n] = k[X_1, \cdots, X_n]$.

  For this purpose, we generalize it to the following form~: 

\noindent
 {\sf  The Generalized Jacobian Conjecture $(GJC)$}~:  
 {\it  Let  $\varphi : S \rightarrow T$ be an unramified homomorphism of Noetherian domains with  $T^\times = \varphi(S^\times)$. Assume that $T$ is a factorial domain and that $S$ is a simply connected normal domain.  Then  $\varphi$ is an isomorphism.   
} 

 For the consistency of our discussion, we raise some serious (or idiot) questions and some comments concerning the examples appeared in the papers published by the certain excellent mathematicians (though we are unwilling to deal with them). Since the existence of such  examples would be against our original target Conjecture$(GJC)$, we have to dispute their arguments about the existence of their respective (so called) counter-examples. Our conclusion is that they are not perfect counter-examples as are shown explicitly in this article.   
 \end{abstract}


\renewcommand{\normalbaselines}{\baselineskip23pt \lineskip4pt \lineskiplimit5pt}
\newcommand{\longmapright}[1]{\smash{\mathop{\hbox to 1.1cm{\rightarrowfill}}\limits^{#1}}}
\newcommand{\mapright}[1]{ \xrightarrow{#1}}
\newcommand{\mapleft}[1]{ \xleftarrow{#1}}
\newcommand{\mapdown}[1]{ \Big\downarrow\llap{ $\vcenter{\hbox{$\scriptstyle#1\,$}}$ } }
\newcommand{\mapup}[1]{ \Big\uparrow\llap{ $\vcenter{\hbox{$\scriptstyle#1\,$}}$ } }

\newcommand{\Liminj}[1]{\raisebox{-1.98ex}%
                        {$\stackrel{\normalsize\mbox{$\varinjlim$}}{\scriptstyle #1}$}}
\newcommand{\Limproj}[1]{\raisebox{-2.00ex}%
                        {$\stackrel{\normalsize\mbox{$\varprojlim$}}{\scriptstyle #1}$}}
             
   
\renewcommand{\thefootnote}{\fnsymbol{footnote}} 

\footnote[0]{
\ \ \  unramified, \'{e}tale, simply connected, factorial and  polynomial rings}

\markboth{\protect\footnotesize{SUSUMU ODA}}{\protect\footnotesize{SOME COMMENTS AROUND THE  GENERALIZED JACOBIAN CONJECTURE}}   

 \bigskip

{\bf Contents}
\begin{small}

\noindent   
\ \ \ 1. Introduction \dotfill \pageref{con2.1}\\
\ \ \ 2. Some Comments around The Generalized Jacobian Conjectures \dotfill \pageref{con04}\\
\ \ \ \ \ \ $\bullet$2.1. Concerning the example of V.Kulikov \dotfill \pageref{con4.1}\\
\ \ \ \ \ \ $\bullet$2.2. Concerning the example of K.Adjamagbo \dotfill \pageref{con4.2}\\
\ \ \ \ \ \ $\bullet$2.3. Concerning the example of F.Oort \dotfill \pageref{con4.3}\\
\ \ \ Appendix A. A Collection of Tools Required in This Paper \dotfill \pageref{con03}\\
\end{small}


\baselineskip=5.32mm 


\vspace{6mm}  
              
{\large \section{Introduction}} \label{con2.1}
  
\bigskip
    
   A fundamental reference for The Jacobian Conjecture $(JC_n)$ is {\bf [6]}.


\bigskip

 Throughout this paper,  unless otherwise specified, we use the following notations~:

\noindent
{\bf $\langle$ Basic Notations and Definitions.$\rangle$}

\noindent
\textbullet\  All fields, rings and algebras are assumed to be  commutative with unity.
 
--- For a ring $ R$,

\noindent
\textbullet\  A {\it factorial domain} $R$ is also called a (Noetherian) unique factorization domain, \\
\textbullet\  $R^{\times}$ denotes the set of units of $R$,\\
\textbullet\  $nil(R)$ denotes the {\it nilradical} of $R$, {\it i.e.,} the set of the nilpotent  elements of $R$,\\
\textbullet\  $K(R)$ denotes the total quotient ring (or the total ring of fractions) of $R$, that is, letting $S$ denote the set of all non-zerodivisors in $R$, $K(R):= S^{-1}R$,\\
\textbullet\  ${\rm Ht}_1(R)$ denotes the set of all prime ideals of height one in $R$,\\
\textbullet\  ${\rm Spec}(R)$ denotes the {\it affine scheme} defined by $R$ (or merely the set of all prime ideals of $R$),\\
\textbullet\  Let $A \rightarrow B$ be a ring-homomorphism and $p\in {\rm Spec}(A)$. Then  $B_p$ means $B\otimes_AA_p$. 

--- Let $k$ be a field.  

\noindent 
\textbullet\   A (separated) scheme over a field $k$ is called {\it $k$-scheme}. A $k$-scheme locally of finite type over $k$ is called   a  {\it $($algebraic$)$ variety} over $k$ or a {\it $($algebraic$)$ $k$-variety} if it is \uline{integral} ({\it i.e.,} irreducible and reduced).\\
\textbullet\   A $k$-variety $V$ is called a {\it $k$-affine variety} or an {\it affine variety over $k$} if it is $k$-isomorphic to an affine scheme ${\rm Spec}(R)$ for some \uline{$k$-affine domain} $R$ ({\it i.e.,} $R$ is a finitely generated domain over $k$).\\
\textbullet\  An integral, closed $k$-subvariety of codimension one in a $k$-variety $V$ is called  a {\it hypersurface} of $V$. \\
\textbullet\   A closed $k$-subscheme (possibly reducible or not reduced) of pure codimension one in  a $k$-variety $V$ is called an (effective) {\it divisor} of $V$, and thus an irreducible and reduced divisor ({\it i.e.,} a prime divisor) is the same as a hypersurface in our terminology.

\bigskip

We consider the following~:
 
\noindent
 {\sf  The Generalized Jacobian Conjecture $(GJC)$}~:  
 {\it  Let  $\varphi : S \rightarrow T$ be an unramified homomorphism of Noetherian domains with  $T^\times = \varphi(S^\times)$. Assume that $T$ is a factorial domain and that $S$ is a simply connected normal domain.  Then  $\varphi$ is an isomorphism.   
} 
\bigskip

 Concerning The Conjecture$(GJC)$, see  {\sf S.ODA~: A Purely Algebraic Short Approach To The Generalized Jacobian Conjecture   (ArXive~:1203.169 v21 [math.AC] -- Nov 2022)}.

\bigskip

Our Main Objective in this paper is to show~:

{\sf `` For the consistency of our discussion about $(GJC)$, we assert that the examples appeared in the papers ({\bf [12], [2]} and {\bf [20]}) published by the certain excellent mathematicians, which  would be against our original target Conjecture$(GJC)$, are imperfect or incomplete counter-examples ''}.

\bigskip

Our general references for unexplained technical terms of Commutative Algebra are {\bf [14]}, {\bf [15]} and {\bf [Eis]}.

\bigskip

  Remark that we often say in this paper that

\begin{quote}

\noindent
{\sf a ring $A$ is called ``simply connected'' if ${\rm Spec}(A)$ is simply connected, and a ring homomorphism $f : A \rightarrow B$ is ``unramified, \'{e}tale, an open immersion, a closed immersion, $\cdots\cdots$" when ``so" is its morphism  ${}^af : {\rm Spec}(B) \rightarrow {\rm Spec}(A)$, respectively.}
\end{quote}


\vspace{6mm}

 \section{Some Comments around the Generalized Jacobian Conjectures}  \label{con04}

\bigskip


 In this section,  we raise some serious (or idiot) questions and some comments concerning the examples appeared in the papers published by the certain excellent mathematicians (though we are unwilling to deal with them). Since the existence of such  examples would be against our original target Conjecture$(GJC)$, we have to dispute their arguments about the existence of their respective (so called) counter-examples. Our conclusion is that they are not perfect counter-examples as are shown explicitly.


\bigskip


The comments treated here influence greatly Conjecture$(GJC)$.
 The author could not accept the core results in {\bf [12]+[2]} and {\bf [20]} (as will be  explained below)  which would be known as counter-examples to  Conjecture$(GJC)$.  
\uline{The discussion here insists that the examples in {\bf [12]+[2]} and {\bf [20]} are not   perfect counter-examples to Conjecture$(GJC)$ definitely.} 
 So he will examine them in details.
 To make sure,  he would like to ask some questions below which should be answered explicitly.

\vspace{9mm}

{\bf $\langle$CONVENTIONS$\rangle$}

In this section, the quoted parts from the applicable papers are written in sanserif letters.   As a (general) rule we quote them from the original papers, but make a few modification on symbolic notations to avoid the confusion (which are clearly understood by their contexts). The underlines and  $(^\ast\ \ )$ in the quoted parts are added by the author. The responsibility for rewritten parts  belongs to him.

\bigskip

Just to make sure, we recall the following~:
In general, a scheme $X$ is called  {\it affine} if $X$ is isomorphic to ${\rm Spec}(A)$ for a ring $A$, and  $X$ is called {\it $k$-affine} or {\it affine over $k$} (where $k$ is a field) if $X$ is isomorphic to ${\rm Spec}(A)$ for an affine ring over $k$ ({\it i.e.,} $A$ is a $k$-algebra of finite type over $k$).


\vspace{6mm}

{\bf \large $\bullet$ 2.1.  Concerning the example in V.S.Kulikov[12]}  \label{con4.1}

\bigskip

 As mentioned in {\bf Introduction}, V.S.Kulikov{\bf [12]} considers the following generalization of The Jacobian Conjecture~:

\bigskip

\noindent  
{\bf Conjecture}(Kul.$GJ_n$)~: {\it Let $X$ be a simply connected 
 algebraic variety over $\mathbb{C}$ and let $F : X \rightarrow \mathbb{A}^n_{\mathbb{C}}$ be a morphism which is \'{e}tale and surjective modulo codimension $2$}\footnote[2]{Let $f : X \rightarrow Y$ be a morphism of schemes. Then $f$ is called  {\it surjective modulo codimension $2$} if the image $f(X)$ intersects every integral closed subscheme of codimension one in  $Y$. When $X = {\rm Spec}(B)$ and $Y ={\rm Spec}(A)$, affine schemes, and $f^*: A \rightarrow B$, then we say that the ring homomorphism  $f^*$ is {\it surjective  modulo codimension $2$}}.
 {\it  Then $F$ is birational.}

\bigskip
He gives a counter-example to {\bf Conjecture} above.  However, he does not show that $X$  is a \uline{$\mathbb{C}$-affine} variety.

\bigskip

{\boldmath$(^\ast \alpha)$} On the page 351 in {\bf [12,\S3]}, V.Kulikov states the following~:
\bigskip

\begin{quote}
{\sf  
Let $D \subseteq \mathbb{C}^n$ be a (possibly reducible) algebraic hypersurface, and let $y$ be a non-singular point of $D$.  Consider   a real plane $\Pi \subseteq \mathbb{C}^n$ intersecting $D$ transversely at $y$.  Let $C \subseteq \Pi$ be a circle of small radius with center at $y$.
 It is well known that the fundamental group $\pi_1(\mathbb{C}^n\setminus D,o)$ is generated by loops $\gamma$ of the following form~:
 $\gamma$ consists of a path $L$ joining the point $o$ with a point $y_1 \in C$, a loop around $y$ along $C$  beginning and ending at $y_1$, and  returning to $o$ along the path $L$ in the opposite direction. Such loops $\gamma$ (and the corresponding elements in $\pi_1(\mathbb{C}^n\setminus D)$) will be called \uline{geometric} \uline{generators}.
}
  
\end{quote}


We use the same notations as in Conjecture(Kul.$GJ_{n}$) basically, but the same notations are used interchangeably in his paper, so that we recite a few notations for our accommodation.  Since  $F : X \rightarrow \mathbb{C}^n$ is  surjective modulo codimension $2$, there exists a hypersurface $D \subseteq \mathbb{C}^n$  such that the restriction $X\setminus F^{-1}(D) \rightarrow \mathbb{C}^n\setminus D$ is a geometric covering of  $d:= {\rm degree}(F)$ sheets (in the $\mathbb{C}$-topology), which is classified by a subgroup $G \subseteq \pi_1(\mathbb{C}^n\setminus D)$ of index $d$.  Lefschetz Theorem implies that $\pi_1(\mathbb{C}^n\setminus D) \cong \pi_1(\mathbb{C}^2\setminus D\cap \mathbb{C}^2)$ by a generic plane section $\mathbb{C}^2 \subseteq \mathbb{C}^n$, and thus Conjecture(Kul.$GJ_{n}$) is equivalent to the following~:
\bigskip

\noindent
 {\bf Conjecture}$(FJ)$: Let $D \subseteq \mathbb{C}^2$ be a $\mathbb{C}$-curve, and let $G$ be a subgroup of finite index in $\pi_1(\mathbb{C}^2\setminus D, o)$ generated by geometric generators.
 Is it true  that $G$ necessarily coincides with  $\pi_1(\mathbb{C}^2\setminus D, o)$~? 
\bigskip

In {\bf [12]}, V.Kulikov gives a {\bf negative} answer to Conjecture(${\rm Kul.}GJ_2$) (and hence $(FJ)$). He constructs  the following counter-example $X \rightarrow \mathbb{C}^2$ to  Conjecture$({\rm Kul.}GJ_2)$\ (See the part {\boldmath$(^\ast  \beta)$} below).   
 However, it  has not been answered whether $X$ is $\mathbb{C}$-affine. (See  Question 1 below.) 
 
\bigskip

{\boldmath$(^\ast \beta)$}  On the pages 355-358 in {\bf [12,\S3]}, V.Kulikov asserts the following~:

\begin{quote}
 {\sf  
{\bf Example :} Let $\overline{D} \subseteq \mathbb{P}^2_{\mathbb{C}}$ denote the curve of degree $4$ with three cusps  defined by $\sum_{i\not= j}X_i^2X_j^2 -2\sum_{i\neq j \neq k \neq i}X_i^2X_jX_k = 0$\  {\rm (cf.{\bf [12,p.358]} and {\bf [9, Chap.4(4.2)]})}, which is given by O.Zariski as the smallest degree curve whose complement has a non-abelian fundamental group {\bf [24,VII \S2]}. Here $\pi_1(\mathbb{P}^2_{\mathbb{C}}\setminus \overline{D})$ is  indeed a non-abelian group of order $12$ generated by geometric generators $g_1, g_2$ satisfying the relations~:
$$ g_1^2 = g_2^2,\ \ g_1^4 = 1,\ \  (g_1g_2)^3 = g_1^2.$$
\uline{Let  $\overline{L} \subseteq \mathbb{P}^2_{\mathbb{C}}$  be a line transversely intersecting  $\overline{D}$ at four points.}  Then there exists a canonical exact sequence~:
$$ 1 \rightarrow K \rightarrow \pi_1(\mathbb{C}^2\setminus D) \rightarrow \pi_1(\mathbb{P}^2_{\mathbb{C}}\setminus \overline{D}) \rightarrow 1,$$
 \uline{where $\mathbb{C}^2 = \mathbb{P}^2_{\mathbb{C}}\setminus \overline{L}$} and $D = \overline{D}\cap \mathbb{C}^2$.
 Here $K \cong \mathbb{Z}$ is central in $\pi_1(\mathbb{C}^2\setminus D)$.
 Let $G$ be a subgroup of  $\pi_1(\mathbb{C}^2\setminus D)$ generated by a pre-image $\overline{g}_1$ of $g_1$, which is of index $3$ and contains $K$.
{\rm  $(^\ast$ Precisely $G=\langle \overline{g}_1 \rangle,  K=\langle \overline{g}_1^4 \rangle \subseteq G$ and $\pi_1(\mathbb{C}^2\setminus D) = G \sqcup \overline{g}_1\overline{g}_2G  \sqcup (\overline{g}_1\overline{g}_2)^2G$. (?))} 
 Then $G$ defines an \'{e}tale morphism $F : X \rightarrow \mathbb{C}^2$, where $X$ is simply connected variety and $F$ has \uline{degree $3$} and is surjective modulo codimension $2$\ (See {\boldmath$(^\ast \gamma)$} below).

We study this example in more detail.  \uline{The morphism  $F$ can be extended to a finite morphism $\widetilde{F} : \widetilde{X} \rightarrow \mathbb{P}^2_{\mathbb{C}}$ of a normal variety $\widetilde{X}$}, and $K\ (^\ast\  \subseteq G \subseteq \pi_1(\mathbb{C}^2\setminus D) = \pi_1(\mathbb{P}^2_{\mathbb{C}}\setminus (\overline{L}\cup \overline{D}))$ \uline{is generated by
 the geometric generator represented by a loop around $\overline{L}$}.
$$\cdots\cdots\cdots\cdots\cdots\cdots\cdots\cdots$$

Therefore there is an exceptional curve $E_1$ of the first kind in  $\widetilde{X}$  and 
\uline{$\widetilde{X}$ is the projective plane $\widetilde{\mathbb{P}}^2_{\mathbb{C}}$ with a point $x = [1:1:1] \in \mathbb{P}^2_{\mathbb{C}}$ blown up}.
$$\cdots\cdots\cdots\cdots\cdots\cdots\cdots\cdots$$
}
\end{quote}

\bigskip

{\boldmath$(^\ast \gamma)$}  In order to justify the above example,  V.S.Kulikov asserts the following on the pages 353-4 of {\bf [12,\S1]}~:

 \begin{quote}
  {\sf
  If $D$ is a divisor of $\mathbb{C}^n\ (n>1)$  and $G \subseteq \pi_1(\mathbb{C}^n\setminus D)$ is a subgroup of finite index $d$ generated by \uline{a part of the set of geometric generators} {\rm  $(^\ast$ of $\pi_1(\mathbb{C}^n\setminus D))$,} then to $G$  there correspond a simply connected algebraic $\mathbb{C}$-variety $X$ and an \'{e}tale morphism $F : X \rightarrow \mathbb{C}^n$ of degree $d$ which is surjective modulo codimension $2$.
 
  In fact, according to [H.Grauert and R.Remmert:Komplex Ra\"{u}me, Math.Ann.136 (1958), 245-318],  to $G$ there correspond a normal variety $\overline{X}$ and a finite morphism   $\overline{F}: \overline{X} \rightarrow \mathbb{C}^n$ such that  $\overline{F}: \overline{X}\setminus \overline{F}^{-1}(D) \rightarrow \mathbb{C}^n\setminus D$ is the unramified covering associated to the inclusion $G \subseteq \pi_1(\mathbb{C}^n\setminus D,o)$.}

{\sf 
 We pick a base point $\overline{o} \in \overline{F}^{-1}(o) \subseteq \overline{X}\setminus \overline{F}^{-1}(D)$, so that $\overline{F}_* : \pi_1(\overline{X}\setminus \overline{F}^{-1}(D),\overline{o}) \rightarrow  G (\subseteq \pi_1(\mathbb{C}^n\setminus D,o))$ is a group-isomorphism.
 $$\cdots \cdots \cdots \cdots \cdots \cdots \cdots \cdots \cdots \cdots$$
For a geometric generator $\gamma \in \overline{F}_*(\pi_1(\overline{X}\setminus\overline{F}^{-1}(D),\overline{o})) = G \subseteq \pi_1(\mathbb{C}^n\setminus D,o)$,  we denote $L_\gamma$ the irreducible component of the divisor $\overline{F}^{-1}(D)$ such that $\overline{\gamma}:=\overline{F}^{-1}_*(\gamma)$ is a loop around $L_\gamma$. Then $L_\gamma$ does not belong to the ramification divisor of the morphism $\overline{F}$. Let $L = \bigcup L_\gamma$ (union over all geometric generators $\gamma \in G$) {\rm $(^\ast$ Note that $L \not= \emptyset$ by the choice of $G$.)}, and let $S$ be a union of the components of the divisor $\overline{F}^{-1}(D)$ not lying in $L$, so that  $\overline{F}^{-1}(D) = L \cup S$.
 Put $X := \overline{X}\setminus S$.    

{\bf Claim.} The variety $X = \overline{X}\setminus S$ is simply connected.

{\it Proof.} The embedding $\overline{X}\setminus \overline{F}^{-1}(D) \subset \overline{X}\setminus S$ induces an  \uline{epimorphism} of groups
$$ \pi_1(\overline{X}\setminus \overline{F}^{-1}(D)) \rightarrow \pi_1(\overline{X}\setminus S) \rightarrow 1.$$
All \uline{geometric generators} from ${\rm \pi}_1(\overline{X}\setminus \overline{F}^{-1}(D))$ lie in the kernel of this epimorphism. On the other hand, $\pi_1(\overline{X}\setminus \overline{F}^{-1}(D))$ is generated by geometric generators. Hence $\pi_1(\overline{X}\setminus S)$ is trivial. $(^\ast$ {\rm So $X$ is simply connected.)}
 $$\cdots \cdots \cdots \cdots \cdots \cdots \cdots \cdots \cdots \cdots$$
}
\end{quote}

He closed his argument by showing $F = \overline{F}|_X : X = \overline{X}\setminus S \rightarrow \mathbb{C}^2$ is \'{e}tale and surjective modulo codimension $2$.  But it is left without proof in {\boldmath$(^\ast \beta)$} and {\boldmath$(^\ast \gamma)$} that $X$ is a $\mathbb{C}$-affine variety.  
  Noting that geometric fundamental groups $\pi_1(\ \ )$  depend only on  topological spaces $(\ \ )$ in the usual $\mathbb{C}$-topology,  by the same reason of asserting ``an epimorphism" above  the group-homomorphism $\pi_1(\overline{X}\setminus S) \rightarrow \pi_1(\overline{X})$ induced from the inclusion $\overline{X}\setminus S\hookrightarrow \overline{X}$ is an epimorphism  and hence $\pi_1(\overline{X})$ is trivial, that is, $\overline{X}$ is simply connected.
  
\bigskip

{\boldmath$(^\ast \delta)$}\ In page 353{\bf [12,\S1]}),
 V.S.Kulikov also asserts the following~:

\begin{quote}
{\sf  
{\bf Lemma.} Let $L$ be a hypersurface in a simply connected variety $V$.  Then $\pi_1(V\setminus L)$ is generated by \uline{geometric generators} (loops aground $L$). 
$$\cdots \cdots \cdots \cdots \cdots \cdots \cdots \cdots \cdots \cdots$$
}
\end{quote}

\bigskip

A (usual) {\it loop} in a topological space $X$ means a continuous map $\omega : [0,1] \rightarrow X$ with $\omega(0) = \omega(1)$. A geometric generator (named by V.S.Kulikov) seems to  mean a loop turning once around an obstruction (See {\boldmath$(^\ast \delta)$}).

 Geometric generators in $\mathbb{C}^n \setminus D$ with an algebraic hypersurface $D$ (possibly reducible) of  $\mathbb{C}^n$ is certainly defined on the page 351 of {\bf [12]}(See {\boldmath$(^\ast \alpha)$}). But we \uline{can not} find a general definition of {\it geometric generators} in $V \setminus H$ or $\pi_1(V\setminus H)$, where  $V$ is  a (even simply connected affine) variety over $\mathbb{C}$ and $H$ is an algebraic hypersurface (possibly reducible) of $V$, which should be given.
 Inferring from the proof of V.S.Kulikov's Lemma above, we can guess its definition, but it is not clear.
  
 So we want to ask the following question~: 
 

 \bigskip

\noindent
{\bf Question 0:} Let $V$ be a (simply connected) normal algebraic $\mathbb{C}$-variety and let $H$ be an effective divisor in $V$. What is the definition of ``a geometric generator" in $V\setminus H$ or $\pi_1(V\setminus H)$~?
 
\bigskip

In  Conjecture(Kul.$GJ_n$), it must be requested that $X$ is $\mathbb{C}$-affine in considering the original Jacobian Conjecture($JC_n)$.

\bigskip

Even if Kulikov's Example is a counter-example to Conjecture$({\rm Kul.}GJ_2)$ as he asserts, it is not shown that it is a counter-example to Conjecture$(GJC)$.
So  we ask the following question~:


\bigskip

\noindent
{\bf Question 1:} In the assertion in {\boldmath$(^\ast\gamma)$} above, is $X = \overline{X}\setminus S$ a \underline{$\mathbb{C}$-affine} variety~? 


\bigskip
Next, referring to {\bf [9,(4.1.2)]}, we can see the following interesting Remark~:
\bigskip

\noindent
{\bf Remark.}  Consider the cone ($\mathbb{C}$-affine surface) 
$$ C:= \{ x \in \mathbb{C}^3\  |\ x_0^3+x_1^3+x_2^3=0 \}$$
and a point $O:=(0,0,0) \in C \subseteq \mathbb{C}^3$.
 It is easy to see that $C$ is contractible to its vertex and has \uline{only one singular point} $O$.  Then $\pi_1(C) = 1$, {\it i.e.,} $C$ is simply connected.  However,  
$$ \pi_1(C\setminus \{ O \}) \twoheadrightarrow H_1(C\setminus \{ O \})\cong \mathbb{Z}^2\oplus \mathbb{Z}/3\mathbb{Z}.$$
 So $C\setminus \{ O \}$ is not simply connected. Note that $C\setminus \{ O \}$ is not a $\mathbb{C}$-affine variety. 
 When the ambient space $C$ is \uline{singular}, the fundamental group $\pi_1(C)=1$ is drastically   changed by removing a point $O$, which is a subset of \uline{codimension 2}.
 If $L$ is a curve on the surface $C$ passing by the point $O$, then the inclusion $C\setminus L \hookrightarrow C\setminus \{ O \}$ induces the surjection $\pi_1(C\setminus L) = \pi_1((C\setminus \{ O \})\setminus (L\setminus \{ O \}) \twoheadrightarrow \pi_1(C\setminus \{ O \}) \not= 1$, which means that $C\setminus L$ is not simply connected. 

\bigskip

 So we ask incidentally the following Questions (a) and (b)~:
 
\bigskip

\noindent
{\bf Question (a):} In general, if $H$ is an effective divisor of a simply connected, non-singular  \underline{$\mathbb{C}$-affine} algebraic variety $V$,  when is $V\setminus H$ simply connected~?  


\bigskip

  Considering  Theorem of van Kampen (which is  seen in a text book on Topology Theory\ ({\it e.g.,}{\bf [9,(4.2.17)]})), an answer to Question (a) above will come out from  the following question~:


\bigskip

\noindent
{\bf Question (b):} Let $H$ be a  hypersurface of a simply connected, non-singular \underline{$\mathbb{C}$-affine} algebraic variety $V$. When can $V\setminus H$ be simply connected~?

\bigskip

\noindent
{\bf Question 2~:} In {\boldmath$(^\ast \gamma)$}, how can the embedding $\overline{X}\setminus \overline{F}^{-1}(D)=\overline{X}\setminus (L\cup S) \hookrightarrow \overline{X}\setminus S (=:X)$ induce an epimorphism of groups $\pi_1(\overline{X}\setminus \overline{F}^{-1}(D)) \rightarrow \pi_1(\overline{X}\setminus S)$~?  Is it trivial according to forgetting some loops   around only $L$~? 
 
\bigskip

Finally we remark that  V.S.Kulikov{\bf [12]} asserts that $\widetilde{X}$ (resp. the surface  $X$) of {\bf Example} in  {\boldmath$(^\ast \beta)$} above is transformed into $\widetilde{\mathbb{P}}^2$, the projective plane with a point $x \in \mathbb{P}^2_{\mathbb{C}}$ blow up (resp. the Kulikov surface $S(R,C,P)$ (named by  K.Adjamagbo) which will be discussed in the next  $\bullet$2.2). 
 K.Adjamagbo[2] shows that $S(R,C,P)$ is $\mathbb{C}$-affine but does not show it is simply connected. So there is a question left~: whether  $S(R,C,P)^\times = \mathbb{C}^\times$ or not (See the next section.).


\vspace{6mm}

{\bf \large $\bullet$ 2.2.  Concerning the example in K.Adjamagbo[2]} \label{con4.2}
 
\bigskip

 Concerning  {\bf Example in $\bullet$1.1{\boldmath$(^\ast \beta)$}} given by V.S.Kulikov{\bf [12]},  K.Adjamagbo {\bf [2]} (even though it is not officially published) informed us the following (under his notations):

\begin{quote}
{\sf 
 Let $P = [1:1:1] \in \mathbb{P}^2 = \mathbb{P}^2_{\mathbb{C}}$, $(X_1,X_2,X_3)$ a system of indeterminates over $\mathbb{C}$,
 $Q_i = 3X^2_i-X_1X_2-X_1X_3-X_2X_3$ for $(1 \leq i \leq 3)$, three quadratic forms defining three conics passing through $P$, 
$\phi$ the morphism from\ $\mathbb{P}^2\setminus \{ P \}$ to\ $\mathbb{P}^2$ whose homogeneous components are defined by the three previous forms $(^*$ it is also regarded as a rational map $\mathbb{P}^2 \dashrightarrow \mathbb{P}^2$ $)$, 
$R$ the Zariski-closure in\ $\mathbb{P}^2$ of the set of ramification points of $\phi$, which is the cubic with a node at $P$ defined by the form $\sum_{i\not=j}X_i^2X_j-6X_1X_2X_3$, 
$Q$ the generic {\rm $(^\ast$ general$)$} linear combination with complex coefficients of the three previous forms $Q_1, Q_2$ and $Q_3$, $C$ the conic of\ $\mathbb{P}^2$ defined by $Q$, passing by $P$ and meeting  \uline{transversely} the cubic $R$ at each point of their intersection\ $(^\ast${\rm including} $P$? {\rm He probably considers the tangent cone at $P$ instead  of the tangent space.}$)$, and such that the image by $\phi$ of the complement of $C$ in\ $\mathbb{P}^2$ is contained in the complement in\ $\mathbb{P}^2$ of a \uline{line} $L$ of\ $\mathbb{P}^2$\ {\rm $(^\ast${\it i.e,} $\phi(\mathbb{P}^2\setminus C) \subseteq \mathbb{P}^2\setminus L = \mathbb{C}^2 \ (P\in C)$ and hence $\overline{\phi^{-1}(L)} \subseteq C)$},  
$\sigma : \widetilde{\mathbb{P}}^2 \rightarrow \mathbb{P}^2$ the blowing-up at the point $P$ of\  $\mathbb{P}^2$, 
$E$ the exceptional curve of\ $\widetilde{\mathbb{P}}^2$\ ({\it i.e.,} $\sigma^{-1}(P)$), 
$\widetilde{R}$ the strict transform of $R$ by $\sigma$, {\it i.e.,} the irreducible curve of\  $\widetilde{\mathbb{P}}^2$ such that $\sigma^{-1}(R) = E \cup \widetilde{R}$, 
$\widetilde{C}$ the the strict transform of $C$ by $\sigma$, {\it i.e.,} the irreducible curve of\  $\widetilde{\mathbb{P}}^2$ such that $\sigma^{-1}(C) = E \cup \widetilde {C}$, 
and $S(R,C,P)$ the complement of $\widetilde{R} \cup \widetilde{C}$ in\  $\widetilde{\mathbb{P}}^2$  {\rm $(^\ast${\it i.e.,} $S(R,C,P)= \widetilde{\mathbb{P}}^2\setminus (\widetilde{R}\cup \widetilde{C}))$}. Then the rational map $\phi\cdot\sigma : \widetilde{\mathbb{P}}^2 \mapright{\sigma} \mathbb{P}^2 \stackrel{\phi}{\dashrightarrow} \mathbb{P}^2$ induces an \uline{unramified} \uline{morphism} $F_{S(R,C,P)} : S(R,C,P) \rightarrow \mathbb{P}^2\setminus L \cong \mathbb{C}^2$, which is of \uline{geometric degree} $3$  and called {\it the Kulikov morphism} on $S(R,C,P)$. 
}
\end{quote}

\bigskip

 V.S.Kulikov asserts that the Kulikov surface $S(R,C,P)$ (which is called such by  K.Adjamagbo) is isomorphic to $X$ in {\bf Example in $\bullet$2.1{\boldmath$(^\ast \beta)$}} above (See {\bf Supplement} below), which is approved by K.Adjamagbo tacitly.  K.Adjamagbo asserts the following properties~:

{\sl
\begin{list}{}{}
\item[{\rm (i)}]  $S(R,C,P)$ is  \uline{$\mathbb{C}$-affine}, non-singular, rational and \uline{factorial}, \item[{\rm (ii)}] $S(R,C,P)$ is \uline{simply connected} $(^*${\rm this depends on V.S.Kulikov{\bf [12]})},
\item[{\rm (iii)}] all its invertible regular functions on $S(R,C,P)$ are constant $(^*$ {\rm its proof depends on  (ii) above)},
\item[{\rm (iv)}] the Kulikov morphism $F_{S(R,C,P)}:S(R,C,P) \rightarrow \mathbb{C}^2$ is  unramified $(^*$ {\rm and hence \'{e}tale)},
\item[{\rm (v)}] $F_{S(R,C,P)}$ is of geometric degree $3$, {\rm $(^\ast$ This means that  $[K(S(R,C,P)):K(\mathbb{C}^2)]=3$)}
\end{list}
}


\bigskip

 The property (i) is seen in K.Adjamagbo {\bf [2]}, explicitly. And the property (iii) is due to (ii) as he explains.


 However, no explicit proof of  the above assertion (ii) seems to appear in  {\bf [2]}, in spite of its concrete construction. 
 
The proof of (ii) should not depend on  V.Kulikov{\bf [12]} (cf. {\bf Example in $\bullet$2.1{\boldmath$(^\ast  \beta)$}} above, which are based on the  surface seen in O.Zariski{\bf [24]}.)  
 unless it is proved that the variety $X$ in {\bf Example in  $\bullet$2.1{\boldmath$(^\ast \beta)$}} is isomorphic to $S(R,C,P)$.

\bigskip

 A more interesting discussion about the geometric fundamental group of the complement of a three-cuspidal quartic curve in $\mathbb{P}^2$  appears in {\bf [9, pp.131-133]}.

\bigskip

\noindent
{\bf Remark 1.}
 We know that there exists a rational map $\psi:\widetilde{\mathbb{P}}^2 \dashrightarrow \mathbb{P}^2$ such that ${\rm Codim}(\widetilde{\mathbb{P}}^2\setminus {\rm dom}(\psi)) \geq 2$ and that $\psi=\phi\cdot \sigma$ as rational maps (See {\bf [Iit,\S1.25} and {\bf \S2.10]}), that is,  $A:=\widetilde{\mathbb{P}}^2\setminus {\rm dom}(\psi)$ is a finite set.   Hence the morphism $\widetilde{\mathbb{P}}^2\setminus E \mapright{\psi} \mathbb{P}^2$ is equal to $\widetilde{\mathbb{P}}^2\setminus E \mapright{\phi\cdot \sigma} \mathbb{P}^2$ and $A \subseteq E$, which means that the rational map $\psi$ is an extension of the rational map $\phi\cdot \sigma$. We may assume that $\psi$ is a rational map such that ${\rm dom}(\psi)$ is maximal among such rational maps.
 {\bf (} By construction, it is easy to see that if $\psi$ is a morphism, the ramification locus of $\psi$ is equal to $\widetilde{R}$ in $\widetilde{\mathbb{P}}^2$.{\bf )}

\bigskip

  {\bf From now on, we use the notation $\psi$ as in {\bf Remark 1}}.

\bigskip

\noindent
{\bf Remark 2.} V.S.Kulikov {\bf [12,\S3]} asserts that $\widetilde{X}=\widetilde{\mathbb{P}}^2$ and  $\widetilde{F} : \widetilde{X} \mapright{\psi} \mathbb{P}^2$ is  a \uline{finite} \uline{morphism} of \uline{$\deg(\psi)=3$} in view of {\bf Remark 1}\ (See also {\bf $\bullet$2.1}{\boldmath$(^\ast \beta)$} mentioned above). 

\bigskip

  To make sure, we would like to ask the following question~:

\bigskip

\noindent
{\bf Question 3:}
 
\begin{list}{}{}
\item[(1)] \uline{How is $S(R,C,P)$ simply connected~?} 
\item[(2)] How is the variety $X$ of {\bf Example in $\bullet$2.1}{\boldmath$(^\ast \beta)$} isomorphic to $S(R,C,P)$~? 
\item[] (If the question (2) is trivial or clearly answered, the question (1) are automatically answered.)
\end{list}

\bigskip

{\bf Author's Assertion about Question 3.}

\bigskip

Now we discuss {\bf Question 3} in more detail to make sure, where we use the same notation freely  as above. In addition, we use the following new notation~: 
\bigskip

\noindent
{\bf Notation} : for a morphism $\theta : A \rightarrow B$ of schemes, $\theta |$ denotes its restriction to a subscheme of $A$.

\bigskip

We shall show that \uline{ $\widetilde{\mathbb{P}}^2\setminus (\widetilde{R}\cup \widetilde{C}) = S(R,C,P)$ is not simply connected}.

 Note that  $\sigma | : \widetilde{\mathbb{P}}^2\setminus E \rightarrow \mathbb{P}^2\setminus \{ P \}$ is an isomorphism.
Hence the inclusion $\widetilde{\mathbb{P}}^2\setminus E \hookrightarrow \widetilde{\mathbb{P}}^2$ induces a surjection $\pi_1(\widetilde{\mathbb{P}}^2\setminus E) \twoheadrightarrow \pi_1(\widetilde{\mathbb{P}}^2)$ and hence $1=\pi_1(\mathbb{P}^2\setminus \{ P \})=\pi_1( \widetilde{\mathbb{P}}^2\setminus E) \twoheadrightarrow \pi_1(\widetilde{\mathbb{P}}^2)$, that is, $\widetilde{\mathbb{P}}^2$ is simply connected. 
 So since $P \in C\cup R$,  we have the isomorphism $\sigma |: \widetilde{\mathbb{P}}^2\setminus (E\cup \widetilde{R}\cup \widetilde{C}) \rightarrow (\mathbb{P}^2\setminus \{ P \})\setminus (R\cup C) = \mathbb{P}^2\setminus (R\cup C)$.  Thus $\sigma |_*: \pi_1(\widetilde{\mathbb{P}}^2\setminus (E\cup \widetilde{R}\cup \widetilde{C})) \rightarrow \pi_1(\mathbb{P}^2\setminus (R\cup C))$ is an isomorphism. 
 The inclusion $\widetilde{\mathbb{P}}^2\setminus (E\cup \widetilde{R}\cup \widetilde{C}) \hookrightarrow \widetilde{\mathbb{P}}^2\setminus (\widetilde{R}\cup \widetilde{C})$ induces a surjection $\pi_1(\widetilde{\mathbb{P}}^2\setminus (E\cup \widetilde{R}\cup \widetilde{C})) \twoheadrightarrow \pi_1(\widetilde{\mathbb{P}}^2\setminus (\widetilde{R}\cup \widetilde{C}))$.    
 Note that $\pi_1(\mathbb{P}^2\setminus (R\cup C)) = \pi_1(\mathbb{P}^2\setminus ((R\cup C)\setminus \{ P \}))$ according to Lemma \ref{D-1}. 
So $\sigma |_{*} : \pi_1(\widetilde{\mathbb{P}}^2\setminus (E\cup \widetilde{R}\cup \widetilde{C})) \twoheadrightarrow \pi_1(\widetilde{\mathbb{P}}^2\setminus (\widetilde{R}\cup \widetilde{C})) \mapright{\sigma |_{*}} \pi_1(\mathbb{P}^2\setminus ((R\cup C)\setminus \{ P \})) = \pi_1(\mathbb{P}^2\setminus (R\cup C))$ is an isomorphism. 
 Therefore we have the following isomorphisms~: 
 $$\pi_1(\widetilde{\mathbb{P}}^2\setminus (\widetilde{R}\cup \widetilde{C})) \cong \pi_1(\widetilde{\mathbb{P}}^2\setminus (E\cup \widetilde{R}\cup \widetilde{C})) \cong  \pi_1(\mathbb{P}^2\setminus (R\cup C)).$$   

 By the way, it is seen  in Lemma \ref{D-2} that $\pi_1(\mathbb{P}^2\setminus (R\cup C)) \twoheadrightarrow \pi_1(\mathbb{P}^2\setminus (R\cup C))/[\pi_1(\mathbb{P}^2\setminus (R\cup C)),\pi_1(\mathbb{P}^2\setminus (R\cup C))] = H_1(\mathbb{P}^2\setminus (R\cup C)) \cong \mathbb{Z}$, an infinite cyclic group  (Consider a case :  $k=2, d_1=\deg(R)=3, d_2=\deg(C)=2$ in Lemma \ref{D-2}), where $[\ \ ,\ \ ]$ denotes a commutator group.  
Therefore $\pi_1(\widetilde{\mathbb{P}}^2\setminus (\widetilde{R}\cup \widetilde{C})) \cong \pi_1(\mathbb{P}^2\setminus (R\cup C)) \twoheadrightarrow \mathbb{Z}$ is not trivial, that is, 
\bigskip

\noindent
{\bf Assertion.} 
$S(R,C,P) = \widetilde{\mathbb{P}}^2\setminus (\widetilde{R}\cup \widetilde{C})$ is \uline{not} simply connected. Thus Kulikov's example $X$ (in {\boldmath$(^\ast \beta)$}) is not isomorphic to Adjamagbo's example $S(R,C,P)$.
  
\bigskip

 Consequently, to make sure, we must ask \uline{how $S(R,C,P)$ is simply connected} like $X$ in {\bf Example in $\bullet$2.1}{\boldmath$(^\ast \beta)$}. It is asserted that $X$ is isomorphic to $S(R,C,P)$ in {\bf [12]} and {\bf [2]}, but we can not find its proof in their articles. 
 In any case,  it could be absolutely hard to answer {\bf Question 3} positively.

\vspace{9mm}

\noindent              
{\bf ------ \large (Supplement) ------}

For understanding Kulikov+Adjamagbo's Example, we discuss $S(R,C,P)$ in more detail.

According to {\bf Remark 2}, the rational map $\psi : \widetilde{\mathbb{P}}^2 \dashrightarrow \mathbb{P}^2_{\mathbb{C}}$ and the morphism  $\psi|=F_{S(R,C,P)}:S(R,C,P)=\widetilde{\mathbb{P}}^2\setminus (\widetilde{R}\cup \widetilde{C}) \rightarrow \mathbb{P}^2_{\mathbb{C}}\setminus L=\mathbb{C}^2$  are indeed exist, and  $F_{S(R,C,P)}$ is a morphism of degree $3$.

Then  $(\psi|=F_{S(R,C,P)}) : \widetilde{\mathbb{P}}^2\setminus (\widetilde{R}\cup \widetilde{C}) = S(R,C,P) \rightarrow \mathbb{P}^2_{\mathbb{C}}\setminus L = \mathbb{C}^2$ is  \'{e}tale and  surjective modulo codimension $2$.


  We put 

$x:=X_1,\ y:=X_2,\ z:=X_3$, where $\mathbb{P}^2_{\mathbb{C}} = Proj\mathbb{C}[X_1,X_2,X_3]$ 

 and

 $u:=3x^2-xy-yz-zx,\ v:=3y^2-xy-yz-zx,\ w:=3z^2-xy-yz-zx$

for convenience sake.

 Let $\mathbb{P}_{x,y,z}^2$ and $\mathbb{P}_{u,v,w}^2$ denote $\mathbb{P}^2$, where $(x,y,z)$, $(u,v,w)$ are systems of indeterminates over $\mathbb{C}$, respectively, that is, $\mathbb{P}^2_{x,y,z} = {\rm Proj}\ \mathbb{C}[x,y,z]$ and $\mathbb{P}^2_{u,v,w} = {\rm Proj}\ \mathbb{C}[u,v,w]$. Put $P:=[1:1:1] \in \mathbb{P}_{x,y,z}^2$.  
The morphism $\phi$ is  defined by $\mathbb{P}_{x,y,z}^2\setminus \{ P \}  \ni [x:y:z] \mapsto [3x^2-xy-yz-zx:3y^2-xy-yz-zx:3z^2-xy-yz-zx]=[u:v:w] \in \mathbb{P}_{u,v,w}^2$ and \uline{the Zariski-closure $R$ of its ramification locus} in $\mathbb{P}_{x,y,z}^2$ is defined by the irreducible polynomial $x^2y+xy^2+x^2z+y^2z+xz^2+yz^2-6xyz=0$.  \uline{The Zariski-closure $B$ of the branch locus} $\phi(R\setminus \{ P \}) \subseteq \mathbb{P}_{u,v,w}^2$ (which corresponds to $\overline{D}$ in {\bf Example of $\bullet$2.1)} is defined by the irreducible polynomial $u^2v^2+v^2w^2+w^2u^2 -2uvw(u+v+w)$, and hence the curve  $B = \overline{\phi(R\setminus \{ P \})}$  is a (irreducible) closed subvariety of degree $4$ with three cuspidal singularities in $\mathbb{P}^2_{u,v,w}$. 

Note that 

\bigskip

$ u^2v^2+v^2w^2+w^2u^2 -2uvw(u+v+w) = (3x^2-xy-yz-zx)^2(3y^2-xy-yz-zx)^2+(3y^2-xy-yz-zx)^2(3z^2-xy-yz-zx)^2+(3z^2-xy-yz-zx)^2(3x^2-xy-yz-zx)^2-2(3x^2-xy-yz-zx)(3y^2-xy-yz-zx)(3z^2-xy-yz-zx)(3x^2-xy-yz-zx+3y^2-xy-yz-zx+3z^2-xy-yz-zx) = 12(xy+xz+yz)(x^2y+xy^2+x^2z-6xyz+y^2z+xz^2+yz^2)^2$.

\bigskip

\noindent
That is to say,  

\bigskip

 $\ \ \ \ \ \ \  u^2v^2+v^2w^2+w^2u^2 -2uvw(u+v+w)$
 
 $\ \ \ \ \  \ \  \ \ \  \ = 12(xy+xz+yz)(x^2y+xy^2+x^2z+y^2z+xz^2+yz^2-6xyz)^2\ \ \ \ (*).$

\bigskip

This means that for the inclusion $\mathbb{C}[u,v,w] \hookrightarrow \mathbb{C}[x,y,z]$ (which induces the rational map $\phi : \mathbb{P}^2_{x,y,z} \dashrightarrow \mathbb{P}^2_{u,v,w}$ and the morphism $\mathbb{P}^2_{x,y,z}\setminus \{ P \} \rightarrow \mathbb{P}^2_{u,v,w}$), we have in $\mathbb{C}[x,y,z]$  

\bigskip

 $\ \ \ (u^2v^2+v^2w^2+w^2u^2 -2uvw(u+v+w))\mathbb{C}[x,y,z]$

 $\ \ \ \ \  \ \ \ \ \ \  \ = (xy+xz+yz)(x^2y+xy^2+x^2z+y^2z+xz^2+yz^2-6xyz)^2\mathbb{C}[x,y,z].$

\bigskip

So we have the Zariski-closure 
$$\overline{\phi^{-1}(B)}=V_+(xy+xz+yz)\cup V_+(x^2y+xy^2+x^2z+y^2z+xz^2+yz^2-6xyz)$$
 in $\mathbb{P}^2_{x,y,z}$ by $(*)$, where $R=V_+(x^2y+xy^2+x^2z+y^2z+xz^2+yz^2-6xyz)$. 
 Note here that both $xy+xz+yz$ and $x^2y+xy^2+x^2z+y^2z+xz^2+yz^2-6xyz$ are irreducible in $\mathbb{C}[x,y,z]$ and that the restriction of $\sigma$ gives  $\sigma | : \widetilde{\mathbb{P}}^2\setminus E \cong \mathbb{P}_{x,y,z}^2\setminus \{ P \}$. 

  Since $P \not\in V_+(xy+yz+zx)$, we have 
\bigskip
$$\sigma^{-1}(V_+(xy+yz+zx))\cap E = \emptyset.$$

\bigskip

Moreover we have easily :
\bigskip

\noindent
{\bf Remark 3.} Due to the blowing up at $P$, both $\widetilde{R}$ and $\widetilde{C}$ are non-singular (irreducible) curves in $\widetilde{\mathbb{P}}^2$. Since $V_+(xy+yz+zx)$ is a non-singular (irreducible) curve in $\mathbb{P}^2_{x,y,z}$, $\psi(\sigma^{-1}(V_+(xy+yz+zx))$ is also irreducible in $\mathbb{P}^2_{u,v,w}$. Hence we have
$$\psi(\sigma^{-1}(V_+(xy+yz+zx))=B$$
 because $B$ is an irreducible closed subset of $\mathbb{P}^2_{u,v,w}$ and  $\psi : \widetilde{\mathbb{P}}^2 \rightarrow \mathbb{P}^2_{u,v,w}$ is a proper morphism. It is easy to see that $\sigma^{-1}(V_+(xy+yz+zx))$ and $\widetilde{R}$ intersect at finitely many points in $\widetilde{\mathbb{P}}^2$. 

\bigskip

 Consider a simultaneous equation~:

$$
\left\{ 
\begin{array}{rcl}
u(x,y,z)&=&3x^2-xy-yz-zx=1\\
v(x,y,z)&=&3y^2-xy-yz-zx=0\\
w(x,y,z)&=&3z^2-xy-yz-zx=0.
\end{array}
\right.
$$

This equation has only a solution $[\pm \frac{1}{\sqrt{3}}:0:0]_{x,y,z}=[1:0:0]_{x,y,z}$ in $\mathbb{P}^2_{x,y,z}$, which means that $\overline{\phi^{-1}([1:0:0]_{u,v,w})}=[1:0:0]_{x,y,z}$ in $\mathbb{P}^2_{x,y,z}$ ({\it e.g.} by MATHEMATICA).

 
Put $T:= [1:0:0]_{x,y,z}$ in $\mathbb{P}^2_{x,y,z}$.  Then $T$ satisfies both the equations~:
$$xy+xz+yz=0\ \ {\rm and }\ \ x^2y+xy^2+x^2z+y^2z+xz^2+yz^2-6xyz=0,$$
which means that $T \in (R\setminus \{ P \})\cap V_+(xy+yz+zx) \subseteq \mathbb{P}^2_{x,y,z}\setminus \{ P \}$.  It is easy to see that the tangent line of $V_+(xy+yz+zx)$ at $T$ and that of $R=V_+(x^2y+xy^2+x^2z+y^2z+xz^2+yz^2-6xyz)$ at $T$ are the same line $V_+(y+z)$ in $\mathbb{P}^2_{x,y,z}$.
 Since $\sigma^{-1}(T)$ is a ramified point as was seen above, let $\widetilde{T}:=\sigma^{-1}(T)$. 

Considering  the similar argument for $T_1:=[0:1:0]_{x,y,z}$ and $T_2:=[0:0:1]_{x,y,z}$ in $\mathbb{P}^2_{x,y,z}$,
 we have  $\phi(T_1)=[0:1:0]_{u,v,w},\ \phi(T_2)=[0:0:1]_{u,v,w}$ and put $\widetilde{T}_1:=\sigma^{-1}(T_1),\ \widetilde{T}_2:=\sigma^{-1}(T_2)$.  Moreover  for $i=1,2$, $T_i \in (R\setminus \{ P \})\cap V_+(xy+yz+zx) \subseteq \mathbb{P}^2_{x,y,z}\setminus \{ P \}$,  and the tangent line of $V_+(xy+yz+zx)$ at $T_1$ (resp. $T_2$) and that of $R$ at $T_1$ (resp. $T_2$) are the same line $V_+(x+z)$ (resp. $V_+(x+y)$) in $\mathbb{P}^2_{x,y,z}$.

Then we have the following commutative diagram~:

$$
\begin{array}{ccccc}
\widetilde{\mathbb{P}}^2&{---}&\stackrel{\psi}{--------}&\dashrightarrow &\mathbb{P}^2_{u,v,w}\\
 \bigcup &   &   &   & ||\\
\widetilde{\mathbb{P}}^2\setminus \{ \widetilde{P}_1, \widetilde{P}_2 \}&\mapright{\sigma|} &\mathbb{P}^2_{x,y,z} & \stackrel{\phi}{\dashrightarrow}& \mathbb{P}^2_{u,v,w} \\
\bigcup &   & \bigcup  &   & ||\\
\widetilde{\mathbb{P}}^2\setminus E&\mapright{\sigma|} &\mathbb{P}^2_{x,y,z}\setminus \{ P \} &\mapright{\phi|} & \mathbb{P}^2_{u,v,w} \\
\bigcup &   & \bigcup  &   & \bigcup\\
\sigma^{-1}(V_+(xy+yz+zx))&\mapright{\sigma|}&V_+(xy+yz+zx)&\mapright{\phi|}&B\\
\bigcup &   & \bigcup  &   & \bigcup\\
\{ \widetilde{T}, \widetilde{T}_1, \widetilde{T}_2 \} &\mapright{\sigma|}& \{ T, T_1, T_2 \} &\mapright{\phi|}& \{ \phi(T), \phi(T_1), \phi(T_2) \},
\end{array}
$$

\bigskip

\noindent
where   
 $\widetilde{C} \mapright{\sigma|} C \stackrel{\phi}{\dashrightarrow} L$
 with 
$\widetilde{T}, \widetilde{T}_1, \widetilde{T}_2 \not\in \widetilde{C},\ T, T_1, T_2 \not\in C$  and $\psi(\widetilde{T})=[1:0:0]_{u,v,w},\ \psi(\widetilde{T}_1)=[0:1:0]_{u,v,w},\ \psi(\widetilde{T}_2)=[0:0:1]_{u,v,w} \not\in L$.


Consider the strict transform of $\overline{\phi^{-1}(B)}=V_+(xy+yz+zx)\cup R$ by $\sigma$, noting that  $P\in R$,\ $P \not\in V_+(xy+yz+zx)$ and that $\sigma^{-1}(V_+(xy+yz+zx))\cap E=\emptyset$ as mentioned above.  
 By {\bf Remark 2}, we see that the ramification locus of $\psi$ is $\widetilde{R}$ in  $\widetilde{\mathbb{P}}^2$.  Then noting that  $\psi(\sigma^{-1}(V_+(xy+yz+zx)))=B$ by {\bf Remark 3}, 
  we have  
$$\psi^{-1}(B) = \sigma^{-1}(V_+(xy+yz+zx))\cup \widetilde{R} \ \ \ \ \ \ \ \ \ \ \ (**)$$

 (or\ \ \  $\psi^*(B) = \sigma^{-1}(V_+(xy+yz+zx)) + 2\widetilde{R}$\ \ \  as divisors in ${\rm Div}(\widetilde{\mathbb{P}}^2)$, where $\psi^* : {\rm Div}(\mathbb{P}^2_{u,v,w}) \rightarrow {\rm Div}(\widetilde{\mathbb{P}}^2)$),
 
and
$$\psi^{-1}(B)\cap E = \widetilde{R}\cap E = \{ \widetilde{P}_1,\ \widetilde{P}_2 \}\ \ \ \ \ \ \ \ (***),$$
where $\widetilde{P}_1,\ \widetilde{P}_2 \not\in \widetilde{C}$ because $C$ intersects $R$ transversely at $P$.


 Put $U:=[1:0:0]_{u,v,w} \in B \subseteq \mathbb{P}^2_{u,v,w}$. 
 Then $\phi^{-1}(U)=\{ T \}$ by the preceding argument.  

 Since $L$ is a general line $V_+(au+bv+cw)$ in $\mathbb{P}^2_{u,v,w}$ with $(a,b,c) \in \mathbb{C}^3$, considering $(a,b,c)$ with $a \not= 0$ we can take $L$ such that  $U \not\in L$.
 Then $\phi^{-1}(U) = T \in (R \cap V_+(xy+yz+zx))\setminus C  \subseteq \mathbb{P}^2_{x,y,z}\setminus C$. Moreover it is easy to check that $P$ is a non-singular point of $C=V_+(au(x,y,z)+bv(x,y,z)+cw(x,y,z)) (\subseteq \mathbb{P}^2_{x,y,z})$ for any $(a,b,c) \not= (0,0,0)$ in $\mathbb{C}^3$.


 In addition, it is clear that both the irreducible curves $\sigma^{-1}(V_+(xy+yz+zx))\setminus \widetilde{C}$ and $\widetilde{R}\setminus \widetilde{C}$ in $\widetilde{\mathbb{P}}^2\setminus \widetilde{C}$ are tangent to the curve $\sigma^{-1}(V_+(y+z))\setminus \widetilde{C}$ (resp. $\sigma^{-1}(V_+(x+z))\setminus \widetilde{C},\ \sigma^{-1}(V_+(x+y))\setminus \widetilde{C}$) at the point $\widetilde{T}$ (resp. $\widetilde{T}_1,\ \widetilde{T}_2)$. 

 We see that  $(E\cap \psi^{-1}(U))\setminus \widetilde{C} \subseteq (E \cap \psi^{-1}(B))\setminus \widetilde{C} = \{ \widetilde{P}_1, \widetilde{P}_2 \}\setminus \widetilde{C} = \emptyset$  by  $(***)$  because $\widetilde{P}_1, \widetilde{P}_2$ are not in $E\cap \widetilde{C}$ (See the first quoted part of {\bf $\bullet$2.2}). 
 Thus $(E\cap \psi^{-1}(U))\setminus \widetilde{C} = \emptyset$.

 Note that the point $\widetilde{T}$ belongs to $\psi^{-1}(U)\setminus \widetilde{C}$ by {\boldmath$(a)$} because $U \not\in L$. Moreover $\psi(U)\setminus \widetilde{C} = \sigma^{-1}\phi^{-1}(U)\setminus \widetilde{C}$.  Indeed, if there exists a point in $(\psi^{-1}(U)\setminus \sigma^{-1}\phi^{-1}(U))\setminus \widetilde{C}$, then it is an unramified point of $\psi$ and is in $E$. But it is impossible because $(E\setminus \widetilde{C})\cap (\psi^{-1}(U)\setminus \widetilde{C}) = \emptyset$ as above. 
 Hence $\psi^{-1}(U)\setminus \widetilde{C}$ consists of only the point $\widetilde{T}$.  It follows that $(\psi^{-1}(U)\setminus \widetilde{C})\cap (\widetilde{\mathbb{P}}^2\setminus (\widetilde{R}\cup \widetilde{C})) = \{ \widetilde{T} \} \cap (\widetilde{\mathbb{P}}^2\setminus (\widetilde{R}\cup \widetilde{C})) = \emptyset$, which means that $\psi|: \widetilde{\mathbb{P}}^2\setminus (\widetilde{R} \cup \widetilde{C}) \rightarrow \mathbb{P}^2_{u,v,w} \setminus L = \mathbb{C}^2$ is not surjective, but may be surjective modulo codimension $2$. 

Put $M:=V_+(xy+yz+zx) \subseteq \mathbb{P}^2_{x,y,z}$, and $\widetilde{M}$ denotes the strict transform of $M$ by $\sigma$, that is, $\widetilde{M} = \overline{\sigma^{-1}(M)\setminus E}$ in $\widetilde{\mathbb{P}}^2$.
 Note that $P\not\in M$ and hence that $E\cap \widetilde{M} = \emptyset $ in $\widetilde{\mathbb{P}}^2$.  So $\widetilde{M}=\sigma^{-1}(M)$ and $\psi(\widetilde{M}) = \phi(M) = B$ by {\bf Remark 3}.

\bigskip

Summarizing the preceding arguments,  we picture the situation roughly as the following figures~:

\begin{small}

\unitlength 0.1in
%

\end{small}

\bigskip

\noindent
where $C\cap M=\{ G_1, H_1, I_1, J_1 \}$,\ $C\cap R=\{ G_2, H_2, I_2, J_2, P \}$,\ 
 $\widetilde{C}\cap \widetilde{M} = \{ \widetilde{G}_1, \widetilde{H}_1, \widetilde{I}_1, \widetilde{J}_1 \}$,\ $\widetilde{C}\cap \widetilde{R}=\{ \widetilde{G}_2, \widetilde{H}_2, \widetilde{I}_2, \widetilde{J}_2 \}$ and $\sigma(\widetilde{P}_1)=\sigma(\widetilde{P}_2)=\sigma(\widetilde{P}_3)=P$ ;

\bigskip


\begin{small}
\unitlength 0.1in
%
\end{small}

\bigskip

\noindent
where $\psi(\widetilde{G}_1)=\psi(\widetilde{G}_2)=G,\ \psi(\widetilde{H}_1)=\psi(\widetilde{H}_2)=H,\ \psi(\widetilde{I}_1)=\psi(\widetilde{I}_2)=I,\ \psi(\widetilde{J}_1)=\psi(\widetilde{J}_2)=J$, \  $\psi(\widetilde{R})=B,\ \psi(\widetilde{C})=L$ and $\psi(\widetilde{P}_1) \in B$,\ $\psi(\widetilde{P}_2) \in B$,\  $\psi(\widetilde{P}_3) \in L$~;   

\bigskip


\begin{small}
\unitlength 0.1in
%
\end{small}

\bigskip


\vspace{6mm}

{\bf \large $\bullet$ 2.3. Concerning the example in F.Oort[20]}  \label{con4.3}

\bigskip

 F.Oort{\bf [20,\S1 and \S5]} informed us of the interesting thing concerning a simple cover of $\mathbb{A}^1_{\mathbb{C}}$, which could contain a counter-example to Conjecture$(GJC)$.  
However, his information (proof) will be denied in the end due to its incompleteness.

 In this Section, we assume that 

{\bf `` Every algebraic variety is defined over $\mathbb{C}$ ''}.

\bigskip

He asserts the following statements which are rewritten according to the above assumption  (under almost all of his notations except for ${(\ \ )}^\times$ (resp. $\mathbb{C}$) instead of ${(\ \ )}^*$ (resp. $k$))~:
\bigskip

 In {\bf [20,\S1]},

\begin{quote}

{\sf 
\noindent
{\bf (1.2)} Let $L$ be a ``function field in one variable" over $\mathbb{C}$ ({\it i.e.,}  $\mathbb{C} \subseteq L$ is an extension of fields of finite type, of transcendence degree one). We write $C$ for the (unique) algebraic curve defined over $\mathbb{C}$, complete {\rm $(^\ast$ hence projective)}, absolutely irreducible and \uline{nonsingular} with field of rational functions $K(C) = L$.  Let $\Sigma_L$ be the set of (equivalence classes of non-trivial) discrete valuations  on $L$  which are trivial on $\mathbb{C}$. Note that an element of $\Sigma_L$ corresponds with a point of $C$.

Suppose $S \subseteq \Sigma_L$ is a chosen {\it finite} set of points on $C$, and $C^0 = C\setminus S$. We write 
$$R_S := \bigcap_{v \not\in S}\mathcal{O}_{C,v}.$$
If $\#(S) > 0$, this $C^0$ is an \uline{affine curve} over $\mathbb{C}$  $(^\ast$ {\rm cf.{\bf [H1,II(4.1)]})} with coordinate {\rm  $(^\ast$ $\mathbb{C}$-affine)} ring $K[C^0]=R_S$. 
$$\cdots \cdots \cdots \cdots \cdots \cdots \cdots \cdots \cdots \cdots$$
 }
\end{quote}
  
\begin{quote} 
{\sf 

\ \ \  A (ramified) cover $\varphi : C \rightarrow \mathbb{P}^1_{\mathbb{C}}$\  {\rm $(^\ast$ that is,  a {\it covering} defined in {\bf [F]}$)$} of non-singular projective curves (over $\mathbb{C}$) is called {\it simple} if for every point $T \in \mathbb{P}^1_{\mathbb{C}}$ the number of geometric points of $\varphi^{-1}(\{ T \})$ is at least $\deg(\varphi)-1$ ; {\it i.e.,} the cover's ramification is at most of degree $2$, and two ramification points of $\varphi$ in $C$ do not map to the same point of $\mathbb{P}^1_{\mathbb{C}}$. 

 We denote by $S=S(\varphi)=S_\varphi$ the support of the different of $\varphi$, {\it i.e.,} the set of points in $C$ where $\varphi$ is ramified.  If moreover $P \in \mathbb{P}^1_{\mathbb{C}}$, we write  $S_{\varphi,P}$ for the set of points on $C$ either ramified under $\varphi$ or mapping onto $P$ ( {\it i.e.,} $S_{\varphi,P} = S(\varphi)\cup {\rm Supp}(\varphi^{-1}(\{ P \})$), and $R_{S_\varphi,P}$ for the corresponding coordinate ring $K[C^0]$ of $C^0 = C\setminus S_{\varphi,P}$. 
$$\cdots \cdots \cdots \cdots \cdots \cdots \cdots \cdots \cdots \cdots$$
}
\end{quote}

\bigskip
 In {\bf [20,\S5]},

\begin{quote}

{\sf 

\noindent 
  {\bf 5.1. Theorem:} Suppose given integers $g$ and $d$ with $d>g\geq 2$. Then there exists a (non-singular projective) curve $C$ defined over $\mathbb{C}$  and a morphism  
$$\varphi : C \longrightarrow \mathbb{P}^1_{\mathbb{C}},\ \ \mbox{and}\ P \in \mathbb{P}_{\mathbb{C}}^1$$
  such that :

$\bullet$  $\varphi$ is a simple cover,

$\bullet$ $\deg(\varphi) = d$ and {\rm genus}$(C) = g$, 

$\bullet$ $R^\times_{S_{\varphi, P}} = \mathbb{C}^\times$.


\bigskip

\noindent
{\bf 5.2 Corollary:} For every $d \in \mathbb{Z}_{\geq 3}$, there exist a $\mathbb{C}$-affine curve $C^0$ over $\mathbb{C}$ and a morphism 
 $$\pi : C^0 \longrightarrow \mathbb{A}^1_{\mathbb{C}}$$
 such that :

$\bullet$ $\pi$ is surjective {\rm $(^\ast$ and of degree $d$)},

$\bullet$ $\pi$ is \'{e}tale,

$\bullet$ $ K[C^0]^\times = R_{S_{\varphi, P}}^\times = \mathbb{C}^\times$.

{\rm $(^\ast$ $C^0 := C\setminus S_{\varphi,P}$ is $\mathbb{C}$-affine, $\mathbb{A}^1_{\mathbb{C}} = \mathbb{P}^1_{\mathbb{C}}\setminus \{ P \}$ and $\pi = \varphi |_{C^0}$ in 5.1.Theorem.)} 
}
\end{quote}

\bigskip

In order to show the existence of such $C$ (and  $C^0$) in the section 6 of {\bf [20]}, he uses   `` Theory of Moduli '' which is studied in {\bf [F]}.

 He asserts the following in the proof of 5.1.Theorem~:

\begin{quote}

 {\sf  

\noindent
{\bf 6. The proof of 5.1.Theorem.}

In this section we fix an integer $g$ (the genus of $C$), and integer $d$ (the degree of the morphism  $\varphi$ of complete curves, or the degree of $\pi = \varphi^0$ of affine curves over $\mathbb{C}$),  and we suppose $d>g\geq 2$. We write $w=2g-2+2d$ (the number of ramification points in a simple covering).

\bigskip

\noindent
{\bf 6.1. Some moduli spaces.}

We write :
$$\mathcal{X} \longrightarrow \mathcal{R}\times \mathbb{P}_{\mathbb{C}}^1 \longrightarrow \mathcal{H} \longrightarrow (^\ast\  {\rm Spec}(\mathbb{C}))$$
 for the following moduli spaces and forgetful morphisms :
\bigskip

\noindent
Here $\mathcal{H}$ is the Hurwitz scheme :  points of this correspond with (isomorphic classes   of) simple covers $\varphi : C \rightarrow \mathbb{P}_{\mathbb{C}}^1$, where the genus $g$, the degree of $\varphi$ equals $d$, and hence the number of ramification points in $C$, equal to the number of branch points in $\mathbb{P}_{\mathbb{C}}^1$, and this number is equal to $w=2g-2+2d$\ {\rm $(^\ast$ See also {\bf [F,(8.1)]})}.  The functor of simple covers is representable, {\it i.e.,} $\mathcal{H}$ exists, and it is a fine moduli scheme,  the functor and this scheme denoted by $\mathcal{H}= \mathcal{H}^{d,w}$ in {\bf [F]}.  $(^\ast$ {\rm Note that the $\mathbb{C}$-scheme $\mathcal{H}$ is irreducible \ (See {\bf [20,(6.2)],\ [F,(7.5)]})}$)$.
 We write $[\varphi] \in \mathcal{H}$ for the corresponding point in this Hurwitz scheme. 
 {\rm $(^\ast$ cf.{\bf [F,(1.9)]}.)}
\bigskip

\noindent
We denote by $\mathcal{R}$ the scheme representing the functor of simple covers with the ramification points marked, {\it i.e.,} a point of $\mathcal{R}$ is an isomorphism class of $(\varphi,Q_1,\ldots,Q_w)$, where $[\varphi] \in \mathcal{H}$ and the different of $\varphi$ equals $\sum_i Q_i$,
$$ [(\varphi,Q_1,\ldots,Q_w)] \in \mathcal{R}.$$
This functor is representable. Note : if $1 \leq s < t \leq w$  then $Q_s \not= Q_t$.

\noindent
We denote by $\mathcal{X}$ the scheme representing the functor of simple covers with the ramification points marked and the fiber over a point $P \in \mathbb{P}_{\mathbb{C}}^1$ numbered, {\it i.e.,} 
$$[(\varphi,Q_1,\ldots,Q_w, P_1,\ldots,P_d, P)] \in \mathcal{X},$$
  with 
$$[(\varphi,Q_1,\ldots,Q_w)] \in \mathcal{R},\ \ P \in \mathbb{P}_{\mathbb{C}}^1$$
and
 $$C_P = C \times_{\mathbb{P}_{\mathbb{C}}^1} \{ P \} = \sum_{j=1}^d P_j$$
 (as divisors on $C$, we allow points in this fiber above $P$ to coincide). This functor is representable.  {\rm $(^\ast$Note that $\mathcal{R}$ and $\mathcal{X}$ are irreducible \ (See {\bf [20,(6.2)]}.)}
 The morphisms above are the natural forgetful morphisms.  
 $$ \cdots\cdots\cdots\cdots\cdots\cdots\cdots\cdots $$


\bigskip

\noindent
{\bf 6.3.} We come to the proof of 5.1.Theorem. Note that we suppose $k= \mathbb{C}$, and  $d>g \geq 2$.
 
For $[a,b]:= (a_1,\ldots,a_\omega, b_1,\ldots, b_d)  \in \mathbb{Z}^{w + d}$\ \ $(w = 2g-2+2d)$    
$$\Delta_{[a,b]} := \{\ x \in \mathcal{X}\ |\ \sum_{1\leq i \leq w}a_i\cdot Q_i + \sum_{1\leq j \leq d}b_j\cdot P_j \sim 0\ \} \subseteq \mathcal{X},$$
where $x=[(\varphi, Q_1,\ldots,Q_w,P_1,\ldots,P_d,P)] \in \mathcal{X}$, 
 is \uline{a Zariski-closed subset} of $\mathcal{X}$. 

\bigskip

\noindent
{\bf 6.4. Claim.} If $\Delta_{[a,b]} = \mathcal{X}$, then $[a,b] = (0, \ldots, 0)$.

{\rm $(^\ast$ This claim means that $\Delta_{[a,b]}$ with $[a,b]\neq (0,\ldots, 0)$ is a proper subset of $\mathcal{X}$).}

$$\cdots \cdots \cdots \cdots \cdots \cdots \cdots \cdots \cdots \cdots$$

\bigskip

\noindent
{\bf 6.5.} \ Proof of 5.1.Theorem : We have seen that for every $[a,b]$ with $[a,b] \not= 0 = (0,\ldots, 0)$ the closed subset $\Delta_{[a,b]} \subset \mathcal{X}$ is a proper subset.  Consider the projection $\Delta'_{[a,b]} \subset \mathcal{H}\times \mathbb{P}^1_{\mathbb{C}}$.  It follows that every $\Delta'_{[a,b]}$  is a proper closed subset of $\mathcal{H}\times \mathbb{P}_{\mathbb{C}}^1$\ (equivalence does not depend  on the order of summation of the points).  Hence
$$\bigcup_{[a,b]\not= 0} \Delta'_{[a,b]}(\mathbb{C}) \not= (\mathcal{H}\times \mathbb{P}^1_{\mathbb{C}})(\mathbb{C}), $$
here we use that a countable union of proper, closed subsets over an uncountable field is still a proper subset.  Choose $(\varphi,P) \in (\mathcal{H}\times \mathbb{P}^1_{\mathbb{C}}(\mathbb{C})$  out side all $\Delta'_{[a,b]}(\mathbb{C})$, with $[a,b]\not= 0$.  Then the curve $C^0=C\setminus S_{\varphi, P}$ has no non-trivial units : every non-trivial unit $f$ would give a non-trivial principal divisor $(f)$, which expresses a linear equivalence of this divisor with support in $S_{\varphi, P}$, a contradiction with $(\varphi, P) \not\in \Delta'_{[a,b]}(\mathbb{C})$ for every $[a,b]\not= 0$. This proves 5.1.Theorem.
} 
\end{quote}
 
\bigskip
     
 It is well-known that ${\rm Pic}(\mathbb{P}_{\mathbb{C}}^1) \cong \mathbb{Z}$, a free additive group of rank one.

 
 \bigskip

Here we have a question~:
\bigskip

\noindent
{\bf Question 4:} Though $\Delta_{[a,b]}\ ([a,b] \in \mathbb{Z}^{w+d})$ is indeed  a subset in $\mathcal{X}$, is  $\Delta_{[a,b]}$ a \uline{Zariski-closed} subset of $\mathcal{X}$~?  If so, how it a Zariski-closed subset of $\mathcal{X}$~? Is it a trivial fact~? 

(This is a core of his argument!)
 
\bigskip


 It is hoped that more ``concrete'' or ``explicit" explanations, or some more ``explicit" references  will be given.

\bigskip


{\bf Author's Assertion about Question 4.}


\bigskip

First of all, to make sure,  we recall the known fact concerning sheaves of abelian groups on topological spaces, which are seen in  {\bf [T]}, {\bf [I,Chap.II]}, and {\bf [K]}.

\bigskip

\noindent
{\bf Lemma} (cf.{\bf [T,\S1]}, {\bf[I,Chap.II]} and {\bf [K,(4.1)]}).   
  \begin{it} Let $\sigma$ and $\tau$ be sections of a sheaf $\mathcal{F}$ of abelian groups on a topological space $X$. Then the set $\{ x \in X\ |\  \sigma_x = \tau_x \}$ is open in $X$.    
\end{it}

\bigskip

Let $\mathcal{F}$ be a presheaf of abelian groups on a topological space $X$. We introduce the presheaf $D(\mathcal{F})$ of ``discontinuous sections'' of $\mathcal{F}$\ (See {\bf [K,(4.1)]}) as follows~: 

For any open subset $U$ of $X$, define 
$$ D(\mathcal{F})(U) := \prod_{x \in U}\mathcal{F}_x.$$
Thus, a section $\tau$ of $D(\mathcal{F})$ over $U$ is a collection $(\sigma_x)_{x\in U}$, where $\sigma_x$ is an element of the stalk $\mathcal{F}_x$.  The restriction ${\rm res}_V^U : D(\mathcal{F})(U) \rightarrow D(\mathcal{F})(V)$ for an open subset $V$ of $U$ sends $(\sigma_x)_{x\in U}$ to $(\sigma_x)_{x\in V}$, (that is, a projection $\prod_{x\in U}\mathcal{F}_x \rightarrow \prod_{x\in V}\mathcal{F}_x$.) Clearly $D(\mathcal{F})$ is a presheaf on $X$. Moreover it is a sheaf on $X$ by definition of sheaves.


\bigskip

{\bf (Back to our argument about [20].)}
 
Take $x \in \mathcal{X}$.  Then $x=[(\varphi_x, Q_{x 1},\ldots, Q_{x w},P_{x 1},\ldots,P_{x d}, P_x)]$, where $\varphi_x : C_x \rightarrow \mathbb{P}^1_{\mathbb{C}}$ is a simple cover  with  a non-singular curve $C_x$, the ramification locus $\{ Q_{x1},\ldots, Q_{x w} \}$ of $\varphi_x$  and $\phi_x{}^{-1}(P_x) = \{ P_{x 1},\ldots,P_{x d} \}$. 

For an open subset $U$ of $\mathcal{X}$,
 let
$$\mathcal{G}(U) := \prod_{x \in U}{\rm Pic}(C_x).$$
 Then it is clear that $\mathcal{G}$ is a presheaf of abelian (additive) groups on $\mathcal{X}$ with the obvious restriction maps. 

 For any open subset $U$ of $\mathcal{X}$, a section $\sigma \in \mathcal{G}(U)$ is a collection $(\sigma_x)_{x \in U} \in \prod_{x\in U}{\rm Pic}(C_x)$, where $\sigma_x$ is an element of an abelian (additive) group ${\rm Pic}(C_x)$, which yields that $\mathcal{G}_x = {\rm Pic}(C_x)$  for each $x \in \mathcal{X}$. 
 Thus  $\mathcal{G}$ is the same as the presheaf $D(\mathcal{G})$ of `` discontinuous sections'' of $\mathcal{G}$. So $\mathcal{G}$ is a sheaf on $\mathcal{X}$, and consequently is a mono-presheaf of abelian (additive) groups on  $\mathcal{X}$ as mentioned above.
   
Here for a divisor $A \in {\rm Div}(C_x)$,  $[A]$ denotes its class in ${\rm Pic}(C_x)$. 

 For $[a,b]=(a_1,\ldots, a_w, b_1,\ldots,b_d) \in \mathbb{Z}^{w+d}$ and for $x \in \mathcal{X}$, let
 $$[[a,b]]_x := (\sum_ia_i\cdot [Q_{x i}]+\sum_jb_j\cdot [P_{x j}]) \in {\rm Pic}(C_x).$$
Then  for an open subset $U$ of $\mathcal{X}$, 
 $$[[a,b]]|_U := ([[a,b]]_x)_{x\in U} \in \mathcal{G}(U)$$
 is a section of $\mathcal{G}$ over $U$, where  $[[a,b]] \in \mathcal{G}(\mathcal{X})$. 

Thus  by {\bf Lemma} above, we  have
\bigskip

\noindent
{\bf Assertion.}
$ \Delta_{[a,b]} = \{ x \in \mathcal{X}\ |\ [[a,b]]_x = 0_x\ {\rm in}\ {\rm Pic}(C_x)  \}$
 is an \uline{Zariski-open} subset of $\mathcal{X}$.
\bigskip

  Therefore for any $[a,b] \in \mathbb{Z}^{w+d}\setminus \{ 0 \}$,  $\Delta_{[a,b]}$ is \underline{Zariski-open} in the irreducible $\mathbb{C}$-scheme $\mathcal{X}$.
 So contrary to F.Oort's argument {\bf (6.3)-(6.5)} above, it can not  necessarily be approved that $\mathcal{X} \supsetneq \bigcup \Delta_{[a,b]},\ [a,b] \in \mathbb{Z}^{w+d}\setminus \{ 0 \}$.  
 So the proofs of {\bf 5.1.Theorem} and {\bf 5.2.Corollary} in {\bf [20]} are incomplete.


\bigskip


\vspace{6mm}

\appendix  

  {\large \section{\bf A Collection of Tools Required in This Paper}}  \label{con03}

\bigskip
Recall the following well-known results, which are required in this paper.  We write down them for convenience.

\bigskip


\begin{lemma}[{\bf [9,Prop(4.1.1)]}] \label{D-1}
 Let $W$ be a  (possibly, reducible) quasi-projective subvariety of $\mathbb{P}^n_{\mathbb{C}}$ and let $\overline{W}$ be its closure. Then the following hold~:\\
{\rm (i)} $\pi_1(\mathbb{P}^n_{\mathbb{C}}\setminus W) = 0$ if $\dim(W)<n-1$~;\\
{\rm (ii)} $\pi_1(\mathbb{P}^n_{\mathbb{C}}\setminus W) = \pi_1(\mathbb{P}^n_{\mathbb{C}}\setminus \overline{W})$ if $\dim(W)=n-1$.
\end{lemma}

\begin{lemma}[{\bf [9,Prop(4.1.3)]}] \label{D-2}
 Let $V_i\ (1 \leq i \leq k)$ be different hypersurfaces of $\mathbb{P}^n_{\mathbb{C}}$ which have   $\deg(V_i)=d_i$. Let $V:= \bigcup_{i=1}^kV_i$.   Then 
 $$\pi_1(\mathbb{P}^n_{\mathbb{C}}\setminus V)/[\pi_1(\mathbb{P}^n_{\mathbb{C}}\setminus V), \pi_1(\mathbb{P}^n_{\mathbb{C}}\setminus V)] = H_1(\mathbb{P}^n_{\mathbb{C}}\setminus V) = \mathbb{Z}^{k-1}\oplus (\mathbb{Z}/(d_1,\ldots, d_k)\mathbb{Z}),$$
where $(d_1,\ldots, d_k)$ denotes the greatest common divisor and $[\  , \ ]$ denotes a  commutator group.  
\end{lemma}

\begin{corollary} \label{DD-3}
  Let $V_i\ (1 \leq i \leq k)$ be different hypersurfaces of $\mathbb{P}^n_{\mathbb{C}}$ which have   $\deg(V_i)=d_i$. Let $V:= \bigcup_{i=1}^kV_i$. 
 Then   $\mathbb{P}^n_{\mathbb{C}}\setminus V$ is simply connected $\Longleftrightarrow$ $V$ is a hyperplane in $\mathbb{P}^n_{\mathbb{C}}$ $\Longleftrightarrow$ $\mathbb{P}^n_{\mathbb{C}}\setminus V \cong \mathbb{A}^n_{\mathbb{C}}$.
\end{corollary}

\begin{proof} 
 By Lemma \ref{D-2}, $\mathbb{P}^n_{\mathbb{C}}\setminus V$ is simply connected if and only if  $k=1$ and $d_1=\deg(V)=1$ if and only if  $V$ is a hyperplane in $\mathbb{P}^n_{\mathbb{C}}$  if and only if $\mathbb{P}^n_{\mathbb{C}}\setminus V \cong \mathbb{A}^n_{\mathbb{C}}$.
 \end{proof}


 \vspace{6mm}



\end{document}